\documentclass[12pt]{amsart} 

\usepackage{amssymb}
\usepackage[pdftex]{graphicx}

\newtheorem{thmA}{Theorem}

\newtheorem{theorem}{Theorem}[section]
\newtheorem{thm}[theorem]{Theorem}
\newtheorem{lemma}[theorem]{Lemma}
\newtheorem{cor}[theorem]{Corollary}
\newtheorem{prop}[theorem]{Proposition}

\theoremstyle{definition}

\theoremstyle{remark}

\newtheorem{remarks}[theorem]{Remarks}

\numberwithin{equation}{section}

\def\E{\Lambda}

\def\<{\langle}
\def\>{\rangle}

\def\rep{{\rm{Rep}}_A}

\def\-{\overline}

\def\X{\mathcal X}

\def\-{\underline}
\def\Z{\mathbb Z}
\def\N{\mathbb N}
\def\G{\Gamma}
\def\S{\Sigma}

\catcode`\@=11
\def\serieslogo@{\relax}
\def\@setcopyright{\relax}
\catcode`\@=12

\begin{document}

\title{The Schur multiplier, profinite completions and decidability}
 
\author[Martin R. Bridson]{Martin R.~Bridson}
\address{Mathematical Institute, 24-29 St Giles', Oxford OX1 3LB, UK} 
\email{bridson@maths.ox.ac.uk} 

\date{7 April 2008; 28 July 2009}



\thanks{This research was supported
by a Senior Fellowship from the EPSRC of Great Britain.}

\begin{abstract} We fix a finitely presented group $Q$ and
consider short exact sequences $1\to N\to \G\to Q\to 1$ with
$\Gamma$ finitely generated.
The inclusion $N\hookrightarrow\G$ induces a morphism
of profinite completions $\hat N\to \hat \G$. We prove
that this is an  isomorphism for all $N$ and $\G$ if and only
if $Q$ is super-perfect and has no proper subgroups of finite index.

We prove that there is no algorithm that, given a 
finitely presented,  residually finite group $\G$ and a finitely presentable
subgroup $P\hookrightarrow\G$, can determine whether or
not $\hat P\to\hat \G$ is an isomorphism.
\end{abstract}

\maketitle
\section*{Introduction}
The profinite completion of a group $\G$ is the inverse system of the directed system of  finite quotients
of $\G$; it is denoted $\hat\G$.
The natural
map $\G\to \hat \G$ is injective if and only if $\G$ is residually finite.  In \cite{BG}  
Bridson and Grunewald settled a question of Grothendieck \cite{groth} by
constructing pairs of finitely presented, residually finite groups $u: P\hookrightarrow\Gamma$ such
that $\hat u :\hat P\to\hat\Gamma$ is an isomorphism but $P$ is not isomorphic to $\Gamma$. 
Grothendieck's interest in this question was motivated by considerations in algebraic geometry:
he wanted to understand the extent to which the topological fundamental group of a complex
projective variety determines the algebraic fundamental group and vice versa. In \cite{groth} he
proved that if $A\neq 0$ is a commutative ring and $u:\G_1\to\G_2$ is a homomorphism of 
finitely generated groups, then $\hat u:\hat\G_1\to\hat\G_2$ is an isomorphism if and only if
the restriction functor $u^*_A:\rep(\G_2)\to\rep(\G_1)$ is an equivalence of categories, where
$\rep(\G)$ is the category of finitely presented $A$-modules with a $\G$-action.

In this article we address the following related
question: given a short exact sequence of groups
 $$1\to N\overset{\iota} \to G \overset{p}\to Q\to 1,$$ 
with $G$ finitely generated,
what conditions on $Q$ ensure that the induced map of
profinite completions $\hat\iota:\hat N\to\hat G$ is an isomorphism, and
how difficult is
it to determine when these conditions hold? We shall also investigate the
map $\hat P\to \hat G\times \hat G$ induced by the inclusion of the
fibre product $P=\{(g,g') \mid p(g)=p(g')\}$.

It is easy to see that $\hat\iota:\hat N\to \hat G$ is surjective if and only if $Q$ has no non-trivial finite quotients.
The following complementary result plays an important role in \cite{PT}, \cite{BL} and \cite{BG};
for a proof see \cite{BG}, p.366.

\begin{lemma} \label{l:inj}
If $Q$ has no non-trivial finite quotients and $H_2(Q;\Z)=0$, then $\hat\iota : \hat N\to \hat G$ is injective.
\end{lemma}

The K\"unneth formula assures us that if 
$Q$ has no non-trivial finite quotients and $H_2(Q;\Z)=0$, then $Q\times Q$ has the
same properties. Thus, since $P/(N\times N) \cong Q$ and
$(G\times G)/(N\times N)\cong Q\times Q$, the lemma implies that the inclusions
$ P \hookleftarrow N\times N \hookrightarrow G\times G$
induce isomorphisms $\hat P \cong \hat N\times \hat N\cong  \hat G\times 
\hat G$.

\begin{cor}\label{c:FP}
If $Q$ has no non-trivial finite quotients and $H_2(Q;\Z)=0$, then 
$P\hookrightarrow G\times G$ induces an isomorphism of profinite
completions.
\end{cor}

The first purpose of this note is to explain why one cannot dispense with the condition $H_2(Q;\Z)=0$
in these results.

\begin{thmA}\label{t:main}
Let $Q$ be a finitely presented group. In order that $\hat\iota : \hat N\to 
\hat G$
be an isomorphism for all short exact sequences 
$1\to N\overset{\iota}\to G \to Q\to 1$, with $G$ finitely generated, it is necessary and sufficient that
\begin{enumerate}
\item $Q$ has no proper subgroups of finite index, and
\item $H_2(Q;\Z)=0$.
\end{enumerate}
\end{thmA}

In Theorem \ref{t:fprods} we shall prove a similar result with $P\hookrightarrow G\times G$ in place of
 $N\hookrightarrow G$.
These results remain valid if one quantifies over various smaller classes of
short exact sequences. For example, one can restrict to the case where $N$ is finitely generated
and $G$ is both finitely presented and residually finite.

\smallskip

It is an open question of considerable interest to determine if there
exists an algorithm that can determine, given a finite presentation, whether the group 
presented has any non-trivial finite quotients or not.  On the other hand, it is known
that there is no algorithm that, given a finite presentation, can determine
whether or not the second homology of the group presented is trivial \cite{gordon}, \cite{CM}.
By exploiting a carefully crafted instance of this phenomenon we prove:

\begin{thmA}\label{t:decide}
There exists a finite set $\X$ and 
recursive  sequences $(R_n)$ and $(S_n)$ of finite sets of words 
of a fixed finite cardinality in the letters $\X^{\pm 1}$,  
such that:
\begin{enumerate}
\item for all $n\in \N$, the group $\G_n=\<\X\mid R_n\>$ is
residually finite;
\item for all $n\in \N$, the subgroup $P_n\subset\G_n$
generated by the image of $S_n$ is finitely presentable;
\item there is no algorithm that can determine for which
$n$ the map $\hat P_n\to \hat\G_n$ induced by inclusion is
an isomorphism.
\end{enumerate}
\end{thmA}

\section{The proof of theorem \ref{t:main}}

We are considering short exact sequences $1\to N\overset{\iota}\to G\overset{p}\to Q\to 1$
with $G$ finitely generated.

It is clear that $\hat\iota :\hat N\to \hat G$ is onto if and only if $Q$ has no non-trivial finite quotients.
Therefore, in the light of Lemma \ref{l:inj}, the following proposition completes the proof of theorem
\ref{t:main}.

\begin{prop}\label{p:NotInj} Let $Q$ be a finitely presented group that
has no non-trivial finite quotients. If $H_2(Q;\Z)$ is non-trivial, then there is a 
short exact sequence $1\to N\overset{\iota}\to G\overset{p}\to Q\to 1$
of residually finite groups,
with $G$
finitely presented and $N$ finitely generated,
such that $\hat\iota :\hat N\to \hat G$ is not injective.
There is also such a sequence with $G$ finitely generated and free (but
$N$ not finitely generated).
\end{prop}

\begin{proof} We shall construct two short exact sequences,
$1\to N\overset{\iota}\to G\overset{p}\to Q\to 1$ 
and $1\to I \overset{u}\to G\to \tilde Q\to 1$, where
$I$ is contained in $N$ and
$\hat u: \hat I \to \hat G$ is surjective ($\tilde Q$ has no
non-trivial finite quotients). We will show that $N/I$ maps onto some non-trivial finite group,
which implies that the morphism $\hat u_N:\hat I\to \hat N$ induced
by inclusion is not surjective. Since
$\hat u= \hat \iota \circ\hat u_N$ is surjective, it will follow that
$\hat \iota$ is not injective.

The group $\tilde Q$ that we use in this construction is the {\em{universal
central extension}} of $Q$. We remind the reader that
a central extension of a group $Q$ is a group $\tilde Q$ equipped
with a  homomorphism $\pi:\tilde Q\to Q$ whose kernel is central
in $\tilde Q$. Such an extension is universal if, given
any other central extension $\pi':E\to Q$ of $Q$, there is a
unique homomorphism $f:\tilde Q\to E$ such that $\pi'\circ f =\pi$.

The standard reference for universal central extensions is 
 \cite{milnor}, pp. 43-47.
The properties that we need here are these:
$Q$ has a universal central extension $\tilde Q$ if (and only if) $H_1(Q;\Z)=0$;
there is a short exact sequence
$$
1\to H_2(Q;\Z)\to \tilde Q\overset{\pi}\to Q\to 1;
$$
and if $Q$ has no non-trivial finite quotients, then neither does $\tilde Q$. 

Since $Q$ and $H_2(Q;\Z)$ are finitely presented, so is $\tilde Q$.
 We fix a finite presentation $\tilde Q=\<A\mid R\>$ and consider
the associated short exact sequence
$$
1\to I \to F \to \tilde Q\to 1,
$$
where $F$ is the free group on $A$ and
$I$ is the normal closure of $R$.
We fix a finite set $Z\subset F$ whose image in $\tilde Q$ generates 
the kernel of $\pi:\tilde Q\to Q$ and define $N\subset F$ to be the normal
closure of $R\cup Z$. Thus we obtain a short exact sequence
$$
1\to N \to F \to Q\to 1.
$$
By construction, $N/I$ is isomorphic to $H_2(Q;\Z)$. Thus it is a non-trivial finitely generated abelian
group, and hence has a non-trivial finite quotient. 

The proof thus far deals with the case where $G$ is free;
it  remains to prove that we can instead arrange for $G$ to be finitely presented
and $N$ to be finitely generated. For this, one
replaces the short exact sequence $1\to I \to F \to \tilde Q\to 1$ in the above argument
by the short exact sequence $1\to I \to \G \to \tilde Q\to 1$ furnished by
Wise's version of the Rips
construction \cite{rips}, \cite{wise}; cf.~the proof of theorem \ref{t:decide} below.
\end{proof}
 
In the next section we shall need the following variation on theorem \ref{t:main}. The {\em{fibre product}} $P\subset \G\times \G$ associated to
a map $\pi:\G\to Q$ is the subgroup $\{(\gamma,\gamma')\mid \pi(\gamma)=\pi(\gamma')\}$.
In the following proof we shall need the (easy) observation that $P$ is a semi-direct product
$K\rtimes\G^\Delta$, where $K=\ker \pi \times \{1\}$ and $\G^\Delta\subset
\G\times\G$ is the diagonal subgroup; if $K$ is central,
this is a direct product.

\begin{thm}\label{t:fprods}
Let $Q$ be a finitely presented group. In order that, for all short exact sequences 
$1\to N\overset{\iota}\to G \to Q\to 1$, the inclusion of the
associated fibre product induce an isomorphism
$\hat P\to \hat G\times \hat G$,
 it is necessary and sufficient that
\begin{enumerate}
\item $Q$ has no proper subgroups of finite index, and
\item $H_2(Q;\Z)=0$.
\end{enumerate}
\end{thm}

\begin{proof} Corollary \ref{c:FP} establishes the sufficiency of the given conditions.
To establish their necessity we argue, using
 the notation of the preceding proof, that if $Q$ has no non-trivial
finite quotients and $H_2(Q;\Z)\neq 0$, then
the map $\hat I\times\hat I\to \hat P$ induced by
inclusion is not surjective. This is enough because 
 $\hat P\to \hat G\times \hat G$ and $\hat I\times \hat I\to \hat G\times \hat G$
are surjective (cf.~corollary \ref{c:FP}).

$I\times I$ is the kernel
of the map $G\times G\to \tilde Q\times \tilde Q$, which sends $P$ to the
fibre product of $\pi:\tilde Q\to Q$. The semi-direct product decomposition 
described prior to this theorem shows that the fibre product of $\pi$ 
is isomorphic to $H_2(Q;\Z)\times \tilde Q$. Hence $P/(I\times I)$ maps onto a 
non-trivial finite (abelian) group.
\end{proof}

\section{An inability to recognise profinite isomorphisms}

Recall that a group presentation is termed {\em{aspherical}} if the universal
cover of the standard 2-complex of the presentation is contractible. 
In section 3.2 of \cite{mb:karl} an argument of Collins and Miller \cite{CM} is
used to construct a
recursive sequence of finite presentations $\Pi_n'\equiv
\<X\mid \-\Sigma_n\>$ with the following properties; see \cite{mb:karl} lemma 3.4.

\begin{lemma}\label{l:ready} Let  $\E_n$ be the group with
presentation $\Pi_n'\equiv \<X\mid \-\S_n\>$.
\begin{enumerate}
\item For all $n\in\N$, the group $\E_n$ has no non-trivial finite quotients.
\item The cardinality of $\-\S_n$ is independent of $n$, and $|\-\S_n| > |X|$.
\item If $\E_n$ is non-trivial then the presentation $\Pi_n'$ is aspherical.
\item There does not exist an algorithm that can determine for which $n$ the group
$\E_n$ is trivial.
\end{enumerate}
\end{lemma}

\begin{cor}\label{c:h2big} If $\E_n\neq 1$ then  $H_2(\E_n;\Z)\neq 0$.
\end{cor}

\begin{proof} If $\E_n\neq 1$ then the presentation $\Pi_n'$ is aspherical and
hence $H_2(\E_n;\Z)$ may be calculated as the second homology group of the
standard 2-complex
of $\Pi_n'$. This complex has no 3-cells and has
more 2-cells than 1-cells, so its second homology is  a non-trivial free-abelian group.
\end{proof}

The following result is Corollary 3.6 of \cite{mb:karl}; the proof relies
on an argument I learned from Chuck Miller.

\begin{prop}\label{p:tilde} There is an algorithm that, given
a finite presentation $\<A\mid \-B\>$ of a perfect group $H$, will
output a finite presentation $\<A\mid B\>$ for the universal
central extension $\tilde H$ (with $\tilde H\to H$ the map induced
by the identity on $A$). Furthermore, $|B| = |X|(1+|\-B|)$.
\end{prop}
 
\begin{remarks}\label{r:z} (1) The images
in $\tilde H$ of the words $b\in \underline B$ will be central, and together
they generate the kernel of $\tilde H\to H$.

(2) If $H$ has a compact classifying space (as the groups $\E_n$
do), then $\tilde H$ does also; cf.~\cite{mb:karl} proposition 1.3.
\end{remarks}

\noindent{\bf{The Proof of Theorem \ref{t:decide}.}}
Let $(\tilde\Pi_n)$ be the  sequence of finite presentations
obtained by applying the algorithm of Proposition \ref{p:tilde}
to $(\Pi_n')$.

We  follow the main construction of \cite{BG} with the presentations
$\tilde\Pi_n\equiv\<X\mid \S_n\>$ of the groups $\tilde\E_n$
 as input. As in \cite{BG}, Wise's version of the Rips construction
\cite{wise} can be used to construct an algorithm that associates to each finite
group-presentation $\mathcal Q\equiv \< Y\mid \Sigma\>$ a finite presentation 
$\mathcal G \equiv \<Y,a_1,a_2,a_3 \mid \check\Sigma\>$ of a residually finite group $G$ such that
$|\check\Sigma|=|\Sigma| + 6|Y|$. There is
an exact sequence 
$$
1\to I\to G\overset{p}\to Q\to 1,
$$
with $I=\<a_1,a_2,a_3\>$, where $Q$ is the group presented by $\mathcal Q$. 

We use this algorithm to convert the   presentations $\tilde\Pi_n\equiv
\<X\mid \S_n\>$ into
 presentations $\<\X\mid R_n\>$ for the
groups $\Gamma_n:=G_n\times G_n$,
where $\X$ is the disjoint union of two copies of $X\cup\{a_1,a_2,a_3\}$ and $R_n$ is the
obvious set of relations arising from the presentation $G_n = \< X,a_1,
a_2,a_3\mid \check\S_n\>$. Now, 
$|\check\Sigma_n| = |\Sigma_n| + 6|X|$ and $R_n$ consists of two copies of $\check\Sigma_n$
together with $(|X| + 3)^2$ commutators. Thus $|R_n|$ is independent of $n$,  as $|\Sigma_n|$ is
(lemma \ref{l:ready}(2) and proposition \ref{p:tilde}).

Let $N_n$ denote the kernel of the composition $G_n\to \tilde\E_n\to \E_n$.
This is generated by 
$\{a_1,\, a_2,\, a_3\}\cup\{\sigma \mid \sigma\in\-\S_n\}$; see remark
\ref{r:z}(1).

Let $P_n\subset G_n\times G_n$ be the fibre product 
associated to the short exact sequence $1\to N_n\to G_n\to  \E_n\to 1$.
Note that $P_n$ is generated by the set $S_n:=\{(a_1,1),\, (a_2,1),\, (a_3,1)\}\cup\{(\sigma, 1)\mid \sigma\in\-\S_n\}\cup \{(x,x)
\mid x\in X\}$, whose cardinality does not depend on $n$.
Now,  $N_n$ is finitely generated, $G_n$ is finitely presented and
$\tilde\E_n$ has a compact classifying space (remark \ref{r:z}(2)).
The 1-2-3 Theorem of \cite{BBMS} 
states that under these circumstances  $P_n$ is finitely presentable (but it does
not give a finite presentation).

In the proof of theorem \ref{t:fprods}, we showed that $\hat P_n\to
\hat G_n\times\hat G_n$ is not injective if $H_2(\E_n;\Z)\neq 0$, and
corollary \ref{c:h2big} assures us that this is the case if $\E_n\neq 1$. On the other
hand, if $\E_n=1$ then  $P=G_n\times G_n$. 
According to lemma \ref{l:ready}(4),
there is no algorithm that can determine, for given
 $n$, which of these alternatives holds.
Thus theorem \ref{t:decide} is proved.\hfill $\square$

\end{document}